# CONCENTRATION OF PERMANENT ESTIMATORS FOR CERTAIN LARGE MATRICES

By Shmuel Friedland, Brian Rider[1] and Ofer Zeitouni[1]

*University of Illinois, Duke University and University of Minnesota*

Let $A_n = (a_{ij})_{i,j=1}^n$ be an $n \times n$ positive matrix with entries in $[a,b]$, $0 < a \le b$. Let $X_n = (\sqrt{a_{ij}} x_{ij})_{i,j=1}^n$ be a random matrix, where $\{x_{ij}\}$ are i.i.d. $N(0,1)$ random variables. We show that for large $n$, $\det(X_n^T X_n)$ concentrates sharply at the permanent of $A_n$, in the sense that $n^{-1} \log(\det(X_n^T X_n) / \operatorname{per} A_n) \to_{n \to \infty} 0$ in probability.

**1. Introduction.** For a set $F \subset \mathbb{R}$ and integers $n \ge m$, denote by $\mathbf{M}(n,m,F)$ the set of $n \times m$ matrices with entries in $F$. Put $\mathbf{M}(n,F) = \mathbf{M}(n,n,F)$. Let $\mathcal{S}_n$ be the symmetric group of permutations acting on $\{1,\ldots,n\}$. For $A \in \mathbf{M}(n,\mathbb{C})$, the permanent of $A$ is defined as

$$\operatorname{per} A = \sum_{\sigma \in \mathcal{S}_n} a_{1\sigma(1)} a_{2\sigma(2)} \cdots a_{n\sigma(n)}.$$

The permanent of a 0–1 matrix is of fundamental importance in combinatorial counting problems. The computation of the permanent of a 0–1 matrix was shown to be a $\#P$-complete problem [15], and, hence, (under standard complexity-theoretical assumptions) not possible in polynomial time. Since then the focus has shifted to randomized approximation methods. The most fruitful method available at present is that of the Markov chain Monte Carlo. In a recent paper Jerrum, Sinclair and Vigoda [10] refined the Markov chain Monte Carlo method to obtain a fully-polynomial randomized approximation scheme for computing the permanents of arbitrary nonnegative matrix.

A second probabilistic method was derived from the following basic observation. Assume

(1.1) $\{x_{ij}\}$ are independent random variables satisfying
$$E(x_{ij}) = 0, \qquad E(x_{ij}^2) = 1.$$

Received February 2002; revised August 2003.
[1]Supported in part by a grant from the Israel Science Foundation, administered by the Israeli Academy of Sciences.

For $A \in \mathbf{M}(n, \mathbb{R}_+)$, let

(1.2) $$X(A) = (\sqrt{a_{ij}} x_{ij})_{i,j=1}^n, \qquad Z(A) = X(A)^T X(A).$$

Then (see [2]),

(1.3) $$\mathrm{E}(\det Z(A)) = \operatorname{per} A.$$

In other words, $\det Z(A)$ is an unbiased estimator of the permanent. The computational advantage of this estimator lies in the well-known fact that the determinant of a large matrix is fast (polynomial) to compute. If $x_{ij}$ are Bernoulli with $x_{ij} \in \{1, -1\}$, then the above estimator is called the Godsil–Gutman estimator [7]. In [2], Barvinok considers the concentration of the estimator (1.3) in the case $x_{ij}$ are Gaussian, complex Gaussian and quaternionic Gaussian. (Of course, moving from real to complex, quaternion or higher-dimensional Gaussians entails some adjustments in the algorithm's description. Namely, the $x_{ij}^2$ appearing in (1.1) should be replaced with $|x_{ij}|^2$ for an appropriate norm-square, and the determinant which makes up the basic estimator should be redefined accordingly. We refer to [3] for a complete discussion of this point.) More precisely, for any $\delta > 0$, Barvinok shows that

$$\sup_{A_n \in \mathbf{M}(n, [0,b])} \mathrm{P}\left(\frac{1}{n} \log \frac{\det Z(A_n)}{\operatorname{per} A_n} \notin [\log \gamma, \delta]\right) \xrightarrow[n \to \infty]{} 0,$$

where $\gamma \approx 0.28$ if $x_{ij}$ are Gaussian, $\gamma \approx 0.56$ if $x_{ij}$ are complex Gaussian and $\gamma \approx 0.76$ if $x_{ij}$ are quaternionic Gaussian. In a more recent preprint [3], Barvinok suggests the possibility of taking $\log \gamma$ any small negative number if each $x_{ij}$ is replaced by a $k \times k$ random matrix with Gaussian entries provided that $k$ is a large enough integer. Along these lines, the work of [4] chooses $x_{ij}$ to be random signed basis elements of a Clifford algebra (of dimension on the order of $n^2$) and proves that in this case $\mathrm{E}[\det Z(A_n)^2]/\mathrm{E}[\det Z(A_n)]^2$ is bounded independently of $n$. Such control of the second moment of the estimator provides concentration via Chebyshev's inequality. Further, since Clifford algebras have representations in terms of real, complex or quaternion matrices of appropriate size, the results of [4] imply that there is very good concentration for real matrices of dimension polynomial in $n$ when the matrices are selected from a set of basis matrices. However, it remains an open question whether this Clifford algebra estimator can be efficiently computed at large dimension. In a sense, both [3] and [4] are guided by the same principle: introducing more randomness at the level of the entries should produce additional averaging and so sharpen the concentration.

In the present note we take a different approach. Our goal is to show that, in fact, good concentration is already present with $k = 1$, if one is willing to look at a restricted class of matrices. In particular, we consider the case where the entries of $A_n$ are uniformly bounded both from above and away from zero. We do this in a slightly more general framework, considering rectangular as well as square matrices. A consequence of our main result, Theorem 2.1, is the following:

COROLLARY 1.1. *Let $0 < a \leq b$ be given. Assume $\{x_{ij}\}$ are independent identically distributed $N(0,1)$ random variables. Then, for any $\delta > 0$,*

$$\lim_{n \to \infty} \sup_{A_n \in \mathbf{M}(n,[a,b])} \mathrm{P}\left(\frac{1}{n}\left|\log \frac{\det Z(A_n)}{\operatorname{per} A_n}\right| > \delta\right) = 0,$$

*where $Z(A_n)$ is defined by* (1.2).

Note that while the restriction $a > 0$ is stronger than one would like (it precludes the important 0–1 matrices), some sort of condition on the entries of the matrices is needed as the the example of $A_n = I_n$ shows. We will have more to say on this point later on.

It has recently been pointed out to us that for the case considered here, namely, with entries bounded above and below, the algorithm of [11] can be adapted to yield a polynomial time $(O(n^4))$ algorithm with polynomially bounded error for computing the permanent. Still, we believe there is an intrinsic interest in the present analysis of Barvinok's algorithm. On one hand, there is the inherent simplicity of the algorithm, with worst case performance bounds, and our results give improved performance for a restricted class of matrices. On the other hand, a study of the algorithm's performance leads directly to rather delicate questions regarding the spectrum of a certain class of random matrices. Indeed, our proof of the above corollary is based on recent concentration results for *nice* functionals of the spectral measure of random matrices [8]. However, since the function $\log(\cdot)$ is not nice enough (it is not globally Lipschitz), a more detailed analysis has to be performed to evaluate the behavior of the bottom or so-called hard edge of the spectrum of $Z(A_n)$. This analysis, which is inspired by ideas of Bai [1], introduces some refinements of current concentration techniques which, we believe, are interesting in their own right and may be applied in other contexts. Indeed, followers of the random matrix theory literature will recognize the $Z(A_n)$ matrices considered here as a natural class of perturbations of the well-known Wishart or Laguerre ensembles.

The structure of the paper is as follows. In Section 2 we introduce our general model of rectangular matrices, state our main theorem, present the basic concentration result we need, and show how the main theorem follows as soon as an integrability condition of the lower tail of the spectrum of $Z(A_n)$ is verified, see Condition 2.1. Section 3 is devoted to the verification of Condition 2.1 under appropriate assumptions on the entries of $A_n$. In Section 4 we study the *flat case* $J_{nm} := M(n, m, [1, 1])$, $n \geq m$. Of course, in this case $\operatorname{per} J_{nm} = \frac{n!}{(n-m)!}$. Our purpose is to point out that for rectangular $J_{nm}$, where $m \leq \theta n$, $\theta < 1$, one immediately gets better concentration than Corollary 1.1. For general, $n \geq m$, we show that a simple polynomial sampling approximates $\operatorname{per} J_{nm}$ to within an error of order one, also tighter than the result in Corollary 1.1. Finally, in the Appendix, we present a more complete study of the lower tail of the spectrum of $Z(J_n)$ by taking advantage of its integrable structure. This analysis, which possesses independent interest, reveals that our Condition 2.1 needed in the course of the proof of Theorem 2.1 is arguably a mild condition.

**2. Preliminaries and main result.** Let $A \in \mathbf{M}(n, m, \mathbb{R}_+)$ (recall that then $m \leq n$). Let also $Q_{m,n}$ denote the set of all strictly increasing sequences $\alpha = \{\alpha_1, \ldots, \alpha_m\} \subset \langle n \rangle$, where $\langle n \rangle = \{1, 2, \ldots, n\}$ and set $A[\alpha, \langle m \rangle] = (a_{\alpha_i j}) \in \mathbf{M}(m, \mathbb{R}_+)$. Then, we define the *permanent* of $A$ as

$$\operatorname{per} A = \sum_{\alpha \in Q_{m,n}} \operatorname{per} A[\alpha, \langle m \rangle].$$

If $A$ is a 0–1 matrix, then $\operatorname{per} A$ counts the matchings of the corresponding bipartite graph.

For $A \in \mathbf{M}(n, \mathbb{R}_+)$, and random variables $\{x_{ij}\}$ satisfying (1.1), the identity (1.3) is immediate, see, for example, [2]. In fact, (1.3) extends to the rectangular case. Indeed, for $A \in \mathbf{M}(n, m, \mathbb{R}_+)$, define $X(A) = (\sqrt{a_{ij}} x_{ij})$ and $Z(A) = X(A)^T X(A)$ as before. Then, using the Cauchy–Binet formula, one finds

$$\operatorname{E}(\det Z(A)) = \operatorname{E}\left(\sum_{\alpha \in Q_{m,n}} \det X[\alpha, \langle m \rangle]^T X[\alpha, \langle m \rangle]\right)$$
$$= \sum_{\alpha \in Q_{m,n}} \operatorname{E}(\det X[\alpha, \langle m \rangle]^T \det X[\alpha, \langle m \rangle])$$
$$= \sum_{\alpha \in Q_{m,n}} \operatorname{per} A[\alpha, \langle m \rangle] = \operatorname{per} A,$$

proving that (1.3) holds true in this case as well.

Our main theorem can now be stated:

THEOREM 2.1. *Let $0 < a \leq b$ be given. Assume $x_{ij}, i, j = 1, \ldots,$ are independent identically distributed $N(0, 1)$ random variables. Then*

$$(2.1) \lim_{n \to \infty} \sup_{A_{n,m} \in \mathbf{M}(n, m, [a, b])} \operatorname{P}\left(\frac{1}{n} |\log \det Z(A_{n,m}) - \log \operatorname{per} A_{n,m}| > \delta\right) = 0,$$

*for any $\delta > 0$.*

Part of the proof of Theorem 2.1 hinges on tailoring certain concentration of measure results to the present setting. To describe these results we must introduce a variety of notations. First, for $\mathbf{S}(m,\mathbb{R}) \subset \mathbf{M}(m,\mathbb{R})$, the set of (real) symmetric matrices, let, for any $B \in \mathbf{S}(m,\mathbb{R})$, denote by $\lambda_1(B) \leq \lambda_2(B) \leq \cdots \leq \lambda_m(B)$ the eigenvalues of $B$ counted with their multiplicities. Recall the spectral factorization $B = QDQ^T$ with orthogonal $Q$, its adjoint $Q^T$ and $D$ the diagonal matrix of eigenvalues. This allows one to view a real valued function $f$ on $\mathbb{R}$ as a function from $\mathbf{S}(m,\mathbb{R})$ into $\mathbf{S}(m,\mathbb{R})$ via $f(B) = Qf(D)Q^T$, where $f(\Lambda)$ is again diagonal with entries $f(\lambda_1(B)), f(\lambda_2(B))$, and so on. And so, along with trace $B = \sum_{i=1}^m \lambda_i(B)$, we may define

$$\text{trace } f(B) = \sum_{i=1}^m f(\lambda_i(B)), \qquad \det f(B) = \prod_{i=1}^m f(\lambda_i(B)).$$

Next, for $f : \mathbb{R} \mapsto \mathbb{R}$, bring in the Lipschitz norm

$$f_{\mathcal{L}} = \sup_{x<y} \frac{|f(x) - f(y)|}{|x-y|};$$

a function $f$ being referred to as Lipschitz when $f_{\mathcal{L}} < \infty$. Lastly, recall that a measure $\nu$ on $\mathbb{R}$ is said to satisfy the logarithmic Sobolev inequality with constant $c$ if, for any differentiable function $f$,

$$(2.2) \qquad \int_{-\infty}^{\infty} f^2 \log \frac{f^2}{\int f^2 \, d\nu} \, d\nu \leq 2c \int_{-\infty}^{\infty} |f'|^2 \, d\nu.$$

The general concentration result of [8], which makes up the backbone of our proof may now be introduced.

Assume that $X \in \mathbf{M}(n,m,\mathbb{R})$ with all $x_{ij}$ mutually independent with laws satisfying the logarithmic Sobolev inequality with uniformly bounded constant $c$. For $Z = (\frac{1}{\sqrt{n}}X)^T(\frac{1}{\sqrt{n}}X)$ and $f$ Lipschitz, Corollary 1.8(b) of [8] states that

$$(2.3) \quad \mathrm{P}\left(\left|\frac{1}{n}\text{trace } f(Z) - \mathrm{E}\left[\frac{1}{n}\text{trace } f(Z)\right]\right| > \delta \frac{n+m}{n}\right) \leq 2\exp\left[-\frac{\delta^2(n+m)^2}{2cf_{\mathcal{L}}^2}\right]$$

for any $\delta > 0$. For us the individual entries of $X = X(A)$ are Gaussian, which are well known to satisfy (2.2). On the other hand, we would like to apply (2.3) with $f = $ the logarithm which is not Lipschitz. This is circumvented by introducing a cutoff: for fixed $\varepsilon > 0$, define

$$\log^\varepsilon x = \log(x \vee \varepsilon),$$

which you will note is Lipschitz. Along with this we set $\det_\varepsilon(B) = \prod_{i=1}^m (\lambda_i(B) \vee \varepsilon)$. Finally, for $A = (a_{ij}) \in \mathbf{M}(n,m,\mathbb{R}_+)$, we define $\tilde{A} = A/\sqrt{n+m}$ and remark that since

$$\frac{\det Z(A)}{\text{per } A} = \frac{\det Z(\tilde{A})}{\text{per } \tilde{A}},$$

it is enough, when proving Theorem 2.1, to consider matrices $\tilde{A}_{n,m} \in \mathbf{M}(n,m, [a/\sqrt{n+m}, b/\sqrt{n+m}])$.

With that said, the form of the general concentration result (2.3) that we will need is stated next as a lemma, the proof of which is deferred to the end of the section.

LEMMA 2.1. *Under the assumptions of Theorem 2.1, let $\varepsilon \in (0, 8b^2)$ and $0 < s_n$, $n = 1, \ldots$, be a sequence diverging to $\infty$. Then, for any $\delta > 0$,*

$$(2.4) \quad \lim_{n \to \infty} \sup_{A_{n,m} \in \mathbf{M}(n,m,[0,b])} \mathrm{P}\bigg(\frac{1}{s_n}|\log \det_\varepsilon Z(\tilde{A}_{n,m}) - \log \mathrm{E}[\det_\varepsilon Z(\tilde{A}_{n,m})]| > \delta\bigg) = 0.$$

*The statement remains true if $\varepsilon = \varepsilon_n \to 0$ as $n \to \infty$, so long as $s_n \varepsilon_n^2 \to \infty$.*

That is, concentration holds at any rate if the small eigenvalues are ignored by way of the cutoff logarithm. Extending beyond the cutoff requires the following integrability condition alluded to above.

Let $\mathcal{A}_{n,m} \subset \mathbf{M}(n,m,[0,b])$.

CONDITION 2.1. *There exist sequences $\varepsilon_n \to 0$, $s_n \to \infty$, such that $s_n \varepsilon_n^2 \to \infty$ and*

$$(2.5) \quad \limsup_{n \to \infty} \sup_{A_{n,m} \in \mathcal{A}_{n,m}} \mathrm{P}\bigg(\frac{1}{s_n} \sum_{\lambda_i(Z(\tilde{A}_{n,m})) < \varepsilon_n} \log \frac{1}{\lambda_i(Z(\tilde{A}_{n,m}))} > \delta\bigg) = 0.$$

Theorem 2.1 is a direct consequence of the following two propositions.

PROPOSITION 2.1. *Fix $b < \infty$ and assume that $\mathcal{A}_{n,m}$ satisfies Condition 2.1. Then, for any $\delta > 0$,*

$$(2.6) \quad \lim_{n \to \infty} \sup_{A_{n,m} \in \mathcal{A}_{n,m}} \mathrm{P}\bigg(\frac{1}{s_n}|\log \det Z(\tilde{A}_{n,m}) - \log \mathrm{E}[\det Z(\tilde{A}_{n,m})]| > \delta\bigg) = 0.$$

PROPOSITION 2.2. *For $0 < a \leq b < \infty$, the class of matrices $\mathbf{M}(n,m,[a,b])$ satisfies Condition 2.1 with $s_n = n$ and $\varepsilon_n = (\log n)^{-4}$.*

Certainly it is of theoretical interest to extend the result of Proposition 2.2 and, thus, Theorem 2.1 to classes of matrices allowing some number of zero entries. As indicated above, there also exist important applied problems for which such a result would be of great use. In this direction we have only the following observation in the case of what we refer to as strictly rectangular matrices. It was pointed out to us, together with its proof, by Silverstein.

PROPOSITION 2.3. *Consider the class of matrices contained in $\mathbf{M}(n, m, [0, b])$ such that $m = m_n$ satisfies $\lim_{n\to\infty}(m/n) = \theta < 1$. Restrict further to the subset of those matrices such that the maximum number of zero entries in any column is bounded by $\lceil \gamma n \rceil$ for all large $n$ with all other entries contained an interval $[a, b]$ bounded away from zero. If, moreover, $\gamma < 1 - \theta$, then Condition 2.1 is satisfied with $s_n = n$ and any $\varepsilon_n \to 0$.*

We conclude this section with the proofs of Proposition 2.1 and its supporting Lemma 2.1. The proofs of Propositions 2.2 and 2.3 are deferred to Section 3.

PROOF OF PROPOSITION 2.1. Take an element $A_{n,m} \in \mathcal{A}_{n,m}$ and for fixed $\delta > 0$ define the numbers

$$g_n(\varepsilon_n, \delta) := \mathrm{P}\left(\frac{1}{s_n} |\log \det_{\varepsilon_n} Z(\tilde{A}_{n,m}) - \log \mathrm{E}[\det_{\varepsilon_n} Z(\tilde{A}_{n,m})]| > \delta\right)$$

and

$$h_n(\varepsilon_n, \delta) := \mathrm{P}\left(\frac{1}{s_n}(\log \det_{\varepsilon_n} Z(\tilde{A}_{n,m}) - \log \det Z(\tilde{A}_{n,m})) > \delta\right)$$

$$= \mathrm{P}\left(\frac{1}{s_n} \sum_{\lambda_i(Z(\tilde{A}_{n,m})) < \varepsilon_n} \log \frac{\varepsilon_n}{\lambda_i(Z(\tilde{A}_{n,m}))} > \delta\right)$$

appearing in the simple bound

(2.7) $$\mathrm{P}\left(\frac{1}{s_n} |\log \det Z(\tilde{A}_{n,m}) - \log \mathrm{E}[\det_{\varepsilon_n} Z(\tilde{A}_{n,m})]| > 2\delta\right)$$
$$\leq g_n(\varepsilon_n, \delta) + h_n(\varepsilon_n, \delta).$$

Next note that, as long as per $\tilde{A}_{n,m} > 0$, one may apply Chebyshev's inequality to the ratio $\frac{\det Z(\tilde{A}_{n,m})}{\operatorname{per} \tilde{A}_{n,m}}$ to produce

(2.8) $$\mathrm{P}\left(\frac{1}{s_n}(\log \det Z(\tilde{A}_{n,m}) - \log \operatorname{per} \tilde{A}_{n,m}) > 2\delta\right) < e^{-2\delta s_n}$$

for $n = 1, 2, \ldots$. Both here and above we are interested in $s_n \uparrow \infty$ while $\varepsilon_n \downarrow 0$.

Now take a small positive $\varepsilon' < 1/4$ and notice that by Condition 2.1 and Lemma 2.1, there exists a large enough integer $N(\delta, \varepsilon')$ so that

$$\sup_{A_{n,m} \in \mathcal{A}_{n,m}} \{g_n(\varepsilon_n, \delta) + h_n(\varepsilon_n, \delta)\} < \varepsilon' \quad \text{and} \quad e^{-2\delta s_n} < \varepsilon' \qquad \text{for all } n > N(\delta, \varepsilon').$$

Hence, for each $A_{n,m} \in \mathcal{A}_{n,m}$, and each $n > N(\delta, \varepsilon')$, the set of $Z(\tilde{A}_{n,m})$ satisfying both inequalities

$$\frac{1}{s_n} |\log \det Z(\tilde{A}_{n,m}) - \log \mathrm{E}[\det_{\varepsilon_n} Z(\tilde{A}_{n,m})]| \leq 2\delta,$$

$$\frac{1}{s_n}(\log \det Z(\tilde{A}_{n,m}) - \log \operatorname{per} \tilde{A}_{n,m}) \leq 2\delta$$

has probability at least $1 - 2\varepsilon'$. Further, since per $\tilde{A}_{n,m} = \mathrm{E}[\det Z(\tilde{A}_{n,m}] \leq \mathrm{E}[\det_{\varepsilon_n} Z(\tilde{A}_{n,m})]$, it follows that

$$\frac{1}{s_n}|\log \mathrm{per}\, \tilde{A}_{n,m} - \log \mathrm{E}[\det_{\varepsilon_n} Z(\tilde{A}_{n,m})]| \leq 4\delta \qquad \text{for } n > N(\delta, \varepsilon').$$

(Note that we deal here with a deterministic difference, hence, if it is bounded above with a positive probability then it is actually bounded above.) Combining the above inequalities with (2.7) we deduce that

$$\sup_{A_{n,m} \in \mathcal{A}_{n,m}} \mathrm{P}\bigg(\frac{1}{s_n}|\log \det Z(\tilde{A}_{n,m}) - \log \mathrm{per}(\tilde{A}_{n,m})| > 6\delta\bigg) \leq \varepsilon'$$

$$\text{for } n > N(\delta, \varepsilon'),$$

completing the proof of Proposition 2.1. □

PROOF OF LEMMA 2.1. Applying (2.3) with the choice $f = \log^\varepsilon$, we obtain

$$\mathrm{P}\Bigg(\bigg|\frac{1}{m+n}\sum_{i=1}^m \log^\varepsilon \lambda_i(Z(\tilde{A}_{n,m})) - \frac{1}{m+n}\mathrm{E}\bigg[\sum_{i=1}^m \log^\varepsilon \lambda_i(Z(\tilde{A}_{n,m}))\bigg]\bigg| \geq \delta\Bigg)$$

(2.9)
$$\leq 2e^{-(m+n)^2 \varepsilon^2 \delta^2/(8b^2)}.$$

The particular form of the right-hand side rests on the readily checked $(\log^\varepsilon)_\mathcal{L} = 1/\varepsilon$ and the well-known fact that a centered Gaussian distribution has logarithmic Sobolev constant equal to its variance. Next set

$$U = \log \det_\varepsilon Z(\tilde{A}_{n,m}) - \mathrm{E}[\mathrm{trace}\, \log^\varepsilon Z(\tilde{A}_{n,m})],$$

and note that (2.9) yields for any $t > 0$,

$$\mathrm{P}(|U| \geq t) \leq 2e^{-\varepsilon^2 t^2/(8b^2)}.$$

Thus,

$$\mathrm{E}[e^U] \leq \mathrm{E}[e^{|U|}] \leq 1 + \int_0^\infty e^t \mathrm{P}(|U| \geq t)\, dt$$

$$\leq 1 + 2\int_0^\infty e^{t - \varepsilon^2 t^2/(8b^2)}\, dt \leq 1 + 2e^{2b^2/\varepsilon^2}.$$

We conclude, together with Jensen's inequality, that

$$\mathrm{E}[\log \det_\varepsilon Z(\tilde{A}_{n,m})] \leq \log \mathrm{E}[\det_\varepsilon Z(\tilde{A}_{n,m})]$$

$$\leq \mathrm{E}[\log \det_\varepsilon Z(\tilde{A}_{n,m})] + \log(1 + 2e^{2b^2/\varepsilon^2}).$$

This, together with (2.9), yields immediately Lemma 2.1 in the case of fixed $\varepsilon$. But by inspecting the above bound, one sees that even if $\varepsilon = \varepsilon_n \downarrow 0$, the statement holds so long as the condition $\varepsilon_n^2 s_n \to \infty$ is respected. □

REMARK 2.1. It is natural to ask what can be said about the performance of the Godsil–Gutman algorithm under the same conditions on the matrix $A$. That is, whether the Gaussians in our statement can be replaced by $\pm 1$ Bernoullis. Toward that, it is true that other concentration results of [8] provide a statement similar to (2.3) as long as the individual laws of the entries of $X$ are compactly supported. However, as we will see in the next section, the isotropic property of the Gaussian is essential in our proof of Proposition 2.2.

**3. Controlling the small eigenvalues.** We now remove the cutoff introduced in the logarithm necessary to go from the concentration inequality of Lemma 2.1 to our main result. That is, the proof of Proposition 2.2 is carried out. In fact, we prove the following slightly stronger statement.

PROPOSITION 3.1. *For all $\varepsilon$ small enough and all $n > m+3$, it holds that*

$$(3.1) \quad \sup_{A_n \in \mathbf{M}(n,m,[a,b])} \mathrm{E}\left[\frac{1}{n} \sum_{\lambda_i(Z(\tilde{A}_n)) < \varepsilon} \log \frac{1}{\lambda_i(Z(\tilde{A}_n))}\right] \leq \frac{\varepsilon|\log \varepsilon|}{a} \frac{(n+m)m}{n(n-m+1)}.$$

*Further, for any $n \geq m$,*

$$(3.2) \quad \limsup_{n \uparrow \infty} \sup_{A_{n,m} \in \mathbf{M}(n,m,[a,b])} \mathrm{E}\left[\frac{1}{n} \sum_{\lambda_i(Z(\tilde{A}_n)) < \varepsilon_n} \log \frac{1}{\lambda_i(Z(\tilde{A}_n))}\right] = 0$$

*as soon as $\varepsilon_n = (\log n)^{-4}$.*

Indeed, Proposition 2.2 follows from (3.2) by Chebyshev's inequality. Note also that in the strictly rectangular case, $\limsup_{n \to \infty}(m/n) < 1$, (3.1) shows that $\varepsilon_n$ may be taken to go to zero arbitrarily slowly. The proof of Proposition 2.3 uses a variant of (3.1); the details are reported at the end of this section.

The proof of Proposition 3.1 makes essential use of the following simple observations.

LEMMA 3.1. *Let $V$ be an element of $\mathbf{M}(n,m,R)$ with statistically independent entries drawn from continuous distributions. Denote by $v_k$ the $k$th column of $V$ and by $V_k$ the matrix formed by deleting $v_k$ from $V$. Then, $\det(V^T V) \neq 0$ and $\det(V_k^T V_k) \neq 0$, a.s. Further,*

$$(3.3) \quad \begin{aligned} \det(V^T V) &= \det(V_k^T V_k)[v_k^T(I - V_k(V_k^T V_k)^{-1} V_k^T)v_k] \\ &=: \det(V_k^T V_k)[v_k^T P_k v_k], \end{aligned}$$

from which it follows that $[(V^T V)^{-1}]_{kk} = (v_k^T P_k v_k)^{-1}$. $P_k$ is a projection, almost surely onto a subspace of dimension $n - m + 1$, and $P_k$ and $v_k$ are independent.

When $V = X(A)$ for $A \in \mathbf{M}(n, m, [a, b])$, one has that $v_k = D_k x_k$ with $D_k = \mathrm{diag}(\sqrt{a_{k1}}, \ldots, \sqrt{a_{kn}})$, and

$$(3.4) \qquad v_k^T P_k v_k = x_k^T D_k P_k D_k x_k = \sum_{i=1}^{n-m+1} \lambda_{i+m-1}(D_k P_k D_k) \hat{x}_{ik}^2,$$

where the $\{\hat{x}_{ik}\}_{i,k=1}^{n,m}$ are independent standard Gaussians and, for each $k$, $\{\hat{x}_{ik}\}_{i=1}^{n}$ and $\{\lambda_i(D_k P_k D_k)\}_{i=1}^{n}$ are also independent. Furthermore, we have the bound

$$(3.5) \qquad a \sum_{i=1}^{n-m+1} \hat{x}_{ik}^2 \leq \sum_{i=1}^{n-m+1} \lambda_{i+m-1}(D_k P_k D_k) \hat{x}_{ik}^2 \leq b \sum_{i=1}^{n-m+1} \hat{x}_{ik}^2.$$

PROOF. The representation (3.3) is commonly exploited in the type of random matrix estimates required below. See, for example, [1] where it is used repeatedly. To understand it, recall the interpretation of $\det V^T V$ as the square of the volume of the parallelepiped spanned by the column vectors $v_1, \ldots, v_m$. Clearly, this is the same as $\det V_k^T V_k$ times the square of the length of the projection of $v_k$ onto the space orthogonal to the span of the columns of $V_k$, but this is just what (3.3) says. That $P_k$ and $v_k$ are independent is clear from the definitions.

Now in the case of $V = X(A)$, one simply notes that the quadratic form $x_k^T D_k P_k D_k x_k$ may be diagonalized by setting $x_k = Q \hat{x}_k$ with an appropriate orthogonal matrix $Q$. By isotropy, the entries of the vector $\hat{x}_k$ remain independent standard Gaussians. The bound on the eigenvalues follows from considering the Rayleigh quotient: with $y = D_k^{-1} z$,

$$a \frac{z^T P_k z}{z^T z} \leq \frac{y^T D_k P_k D_k y}{y^T y} = \frac{z^T P_k z}{z^T D_k^{-2} z} \leq b \frac{z^T P_k z}{z^T z}.$$

From the min–max theorem, one sees that for all $i$,

$$a \lambda_i(P_k) \leq \lambda_i(D_k P_k D_k) \leq b \lambda_i(P_k).$$

As $P_k$ is a projection onto an $n - m + 1$ dimensional subspace, $\lambda_i(P_k) = 0$ for $i \leq m - 1$ and $\lambda_i(P_k) = 1$ for $i > m - 1$, completing the statement. □

PROOF OF PROPOSITION 3.1. We begin with the rectangular case: $A_n \in \mathbf{M}(n, m, [a, b])$, $n > m + 3$ (for economy of space, we will omit the subscript $m$ from $A_{n,m}$). By the monotonicity of $x \log(1/x)$ for $x \in [0, \frac{1}{e}]$, it follows that for any positive $\varepsilon \leq \frac{1}{e}$,

(3.6)
$$\frac{1}{n} \sum_{\lambda_i(Z(\tilde{A}_n))<\varepsilon} \log \frac{1}{\lambda_i(Z(\tilde{A}_n))}$$
$$\leq \frac{\varepsilon|\log \varepsilon|}{n} \sum_{i=1}^{m} \frac{1}{\lambda_i(Z(\tilde{A}_n))} = \frac{\varepsilon|\log \varepsilon|}{n} \sum_{i=1}^{m} [Z(\tilde{A}_n)^{-1}]_{ii} \equiv M_n.$$

By Lemma 3.1 each $[Z(\tilde{A}_n)^{-1}]_{ii}$ is stochastically bounded as in

$$[Z(\tilde{A}_n)^{-1}]_{ii} \leq a^{-1}(n+m)U_{i,n},$$

where $U_{i,n}$ is distributed as one over a $\chi^2$ random variable with $n-m+1$ degrees of freedom, the mean of which one can compute exactly:

$$\mathrm{E}[U_{i,n}] = \frac{\int_0^\infty r^{(n-m-3)/2} e^{-r/2}\, dr}{\int_0^\infty r^{(n-m-1)/2} e^{-r/2}\, dr} = \frac{1}{n-m+1}.$$

Thus, one finds that for $\varepsilon \leq \frac{1}{e}$,

$$\mathrm{E}[M_n] \leq \frac{\varepsilon|\log \varepsilon|}{a} \frac{(n+m)m}{(n-m+1)n},$$

which explains the bound (3.1).

To complete the proof of the proposition, it is enough to consider $n \leq 2m$. Take $A_n \in \mathbf{M}(n,m,[a,b])$, and denote the columns of $X(\tilde{A}_n)$ by $\bar{x}_1,\ldots,\bar{x}_m$ [$\bar{x}_l = (n+m)^{-1/2} D_l x_l$ in previously used notation]. Recall the following identity for the determinant from Lemma 3.1:

$$\det(Z(\tilde{A}_n)) = \det(Z((\tilde{A}_n)_1))[\bar{x}_1^T P_1 \bar{x}_1].$$

By $(A_n)_1$ we mean the matrix formed by the last $m-1$ columns of $A_n$. The matrix $P_1$ projects onto the $(n-m+1)$-dimensional space orthogonal to the span of the columns of $X((A_n)_1)$; it is independent of the vector $\bar{x}_1$.

The above may be iterated: first applying the identity to $\det(Z((\tilde{A}_n)_1))$ and so on. We take a positive $\theta = \theta_n \ll 1$, and after carrying out this procedure $\lceil n\theta_n \rceil$ times, we write the outcome as follows:

$$\det(Z(\tilde{A}_n)) = \det(Z(\tilde{B}_n)) \prod_{k=1}^{\lceil n\theta_n \rceil} [\bar{x}_k^T P_k \bar{x}_k].$$

Here $B_n \in M(n,\bar{m}_n,[a,b])$ is the matrix formed by the last $\bar{m}_n = m - \lceil n\theta_n \rceil$ columns of $A_n$, and each $P_k$ is an $(n-m+k)$-dimensional projection independent of $x_k$. The above equality is re-expressed as

$$\frac{1}{n} \sum_{\lambda_i(Z(\tilde{A}_n))<\varepsilon_n} \log \frac{1}{\lambda_i(Z(\tilde{A}_n))}$$

$$= \frac{1}{n} \sum_{\lambda_i(Z(\tilde{B}_n))<\varepsilon_n} \log \frac{1}{\lambda_i(Z(\tilde{B}_n))} - \frac{1}{n} \sum_{i=1}^{m-\bar{m}_n} \log(\bar{x}_i^T P_i \bar{x}_i)$$

(3.7)
$$+ \frac{1}{n}\left[\sum_{\lambda_i(Z(\tilde{A}_n))\geq\varepsilon_n} \log \lambda_i(Z(\tilde{A}_n)) - \sum_{\lambda_i(Z(\tilde{B}_n))\geq\varepsilon_n} \log \lambda_i(Z(\tilde{B}_n))\right]$$

$$\equiv I_n + II_n + III_n.$$

The point is that the estimates obtained for the strictly rectangular case may be directly applied to $B_n$ and so to $I_n$. That is, we know from (3.1) that there exists a numerical constant $C_1$ such that

$$\text{(3.8)} \qquad \text{E}[I_n] \leq C_1 \frac{\varepsilon_n |\log \varepsilon_n|}{a\theta_n}$$

for all sufficiently large $n$. The term $II_n$ may also be handled by previous considerations. From Lemma 3.1 it follows that: with $a \leq \gamma_i \leq b$ for all $i$ and $\{\hat{x}_i\}$ independent standard Gaussians,

$$-\text{E}[\log(\bar{x}_k^T P_k \bar{x}_k)] = -\text{E}\left[\log\left(\frac{\gamma_1 \hat{x}_1^2 + \cdots + \gamma_k \hat{x}_k^2}{n+m}\right)\right]$$
$$\leq \log(n+m) - \log a - \text{E}[\log \hat{x}_1^2].$$

The last expectation is certainly finite and so there is a constant $C_2$ (depending on $a$ only) such that

$$\text{(3.9)} \qquad \text{E}[II_n] \leq -\left(\frac{\lceil n\theta_n \rceil}{n}\right) \text{E}\left[\log \frac{a\hat{x}_1^2}{n+m}\right] \leq C_2 \theta_n \log n.$$

As for the last term to be bounded, $III_n$, first note that by the interlacing inequalities for any $l \leq \bar{m}_n$,

$$\lambda_l(Z(\tilde{A}_n)) \leq \lambda_l(Z(\tilde{B}_n)) \leq \lambda_{l+\lceil n\theta_n \rceil}(Z(\tilde{A}_n)).$$

Thus, if $l^*$ is the smallest $l$ such that $\lambda_l(Z(\tilde{B}_n)) \geq \varepsilon_n$, then $\lambda_{l^*-1}(Z(\tilde{A}_n)) < \varepsilon_n$ and $\lambda_{l^*+\lceil n\theta_n \rceil}(Z(\tilde{A}_n)) \geq \varepsilon_n$.

Now for each $l$ such that $\lambda_l(Z(\tilde{A}_n)) \geq \varepsilon_n$, the term containing $\log \lambda_l(Z(\tilde{A}_n))$ is paired with the corresponding object in the $B_n$ sum. The contribution to $III_n$ is

$$\frac{1}{n} \log \frac{\lambda_l(Z(\tilde{A}_n))}{\lambda_l(Z(\tilde{B}_n))} \leq 0.$$

In this manner it is possible that the largest $\lceil n\theta_n \rceil$ of the $\lambda_i(Z(\tilde{A}_n))$'s and the smallest $\lceil n\theta_n \rceil$ of the $\lambda_l(Z(\tilde{B}_n))$'s in $III_n$ remain unpaired. That is,

$$III_n \leq \frac{1}{n} \sum_{i=m-\lceil n\theta_n \rceil}^{m} |\log \lambda_i(Z(\tilde{A}_n))| + \frac{1}{n} \sum_{i=l^*}^{l^*+\lceil n\theta_n \rceil} |\log \lambda_i(Z(\tilde{B}_n))|$$

(3.10)
$$\leq 4\theta_n \max\left(|\log \varepsilon_n|, \log \lambda_n(Z(\tilde{A}_n))\right)$$

for all $n$ large enough. The random variable remaining on the right-hand side is, in turn, controlled by

$$\mathrm{E}[\lambda_n(Z(\tilde{A}_n))] \leq \mathrm{E}[\operatorname{trace} Z(\tilde{A}_n)] = \frac{1}{n+m}\mathrm{E}\left[\sum_{i,j} a_{ij} x_{ij}^2\right] \leq bn,$$

an admittedly crude but sufficient bound.

Lastly, (3.8)–(3.10) are combined to produce

$$\mathrm{E}\left[\frac{1}{n}\sum_{\lambda_i(Z(\tilde{A}_n))<\varepsilon_n} \log \frac{1}{\lambda_i(Z(\tilde{A}_n))}\right]$$

$$\leq C_1 \frac{\varepsilon_n |\log \varepsilon_n|}{\theta_n} + C_2 \theta_n \log n + C_3 \theta_n (|\log \varepsilon_n| + \log n)$$

for all large enough $n$. The proof is then finished by choosing $\theta_n = (\log n)^{-2}$ and $\varepsilon_n = (\log n)^{-4}$. □

PROOF OF PROPOSITION 2.3. For simplicity take $m = m_n = \lceil n\theta \rceil$. Tracing the proof of the bound (3.1), one comes to the inequality: with again $A_n = A_{nm_n}$,

$$\frac{1}{n}\mathrm{E}\left[\sum_{\lambda_k(Z(\tilde{A}_n))<\varepsilon} \log \frac{1}{\lambda_k(Z(\tilde{A}_n))}\right]$$

(3.11)
$$\leq \frac{\varepsilon|\log \varepsilon|}{n} \sum_{k=1}^{\lceil n\theta \rceil} \mathrm{E}[Z(\tilde{A}_n)^{-1}]_{kk}$$

$$= \frac{\varepsilon|\log \varepsilon|(n + \lceil n\theta \rceil)}{n} \sum_{k=1}^{\lceil n\theta \rceil} \mathrm{E}\left[\frac{1}{x_k^T(D_k P_k D_k)x_k}\right].$$

Further bounding above requires controlling the eigenvalues of $D_k P_k D_k$ from below. This was previously accomplished by a Raleigh–Ritz argument (Lemma 3.1). In the case that there are some number of zero entries this needs to be replaced by the more sophisticated inequalities of Fan [6].

Note that with the number of zeros in any column bounded by $\lceil n\gamma \rceil$, $P_k$ still projects onto an $(n - \lceil n\theta \rceil + 1)$-dimensional subspace (a.s.). The problem lies in the zeros on the diagonal of $D_k$.

Now for any invertible nonnegative Hermitian matrix $M_1$ and nonnegative Hermitian $M_2$, Fan [6] gives us that

(3.12)
$$\lambda_{i+j+1}^{M_2} \leq \lambda_{i+1}^{M_1 M_2} \lambda_{j+1}^{M_1^{-1}},$$

in which $\lambda_i^M$ is the $i$th largest eigenvalue of the matrix $M$ (twice). By continuity this inequality still holds when $M_1$ has some zero eigenvalues. It is to be applied in this setting with $M_1 = P_k$ and $M_2 = D_k^2$ (the eigenvalues of $D_k P_k D_k$ and $P_k D_k^2$ being the same). With that, (3.12) reduces to

$$\lambda_{\lceil n\theta \rceil + i}^{D_k^2} \leq \lambda_{i+1}^{D_k P_k D_k}.$$

Therefore, if $D_k$ has $n - \lceil \gamma n \rceil$ eigenvalues larger than $a$ (or the $k$th column of $A_n$ has at least that many entries similarly bounded below), then $D_k P_k D_k$ has at least $(n - \lceil n\gamma \rceil - \lceil n\theta \rceil)$ eigenvalues larger than $a$.

By assumption the above holds for each $k$ and $1 - \gamma - \theta > 0$. Thus, the Gaussian quadratic form in the denominator of (3.11) is stochastically larger than $a$ times a $\chi^2$ random variable of degree at least $(1/2)(1 - \theta - \gamma)n$ for all large enough $n$. The right-hand side of (3.11) is then bounded by a constant (depending on $a, \gamma$ and $\theta$) times $\varepsilon |\log \varepsilon|$ and the statement follows. $\square$

**4. The flat case.** This section is devoted to a study of the flat case $A = J_{nm}$, $n \geq m$. This special case is typically referred to as the Laguerre or Wishart Ensemble in the random matrix theory literature. Of course, per $J_{nm}$ is easily computed, and there is no need for an approximate algorithm. However, we wish to emphasize two points in this simpler setting which suggest our general concentration result for the permanent may not be optimal.

The first point focusses on the strictly rectangular case. We have the following:

PROPOSITION 4.1. *Let $n \geq m_n$, $n, m_n \in \mathbb{N}$ and assume that $\{x_{ij}\}_{i,j=1}^{n,m_n}$ are independent identically distributed $N(0,1)$ random variables. Suppose that*

$$\limsup_{n \to \infty} \frac{m_n}{n} \leq \theta < 1. \tag{4.1}$$

*Then for any sequence $s_n$ diverging to $\infty$,*

$$\lim_{n \to \infty} \mathrm{P}\left(\frac{1}{s_n} |\log \det Z(J_{nm_n}) - \log \mathrm{per}\, J_{nm_n}| > \delta\right) = 0.$$

On the other hand, for flat matrices of more general shape we introduce a new polynomial-time estimator that approximates the permanent to within order one error. The statement follows.

PROPOSITION 4.2. *Define $Y_n = n^{-(2+\rho)} \sum_{k=0}^{n^{2+\rho}} X_k^n$, in which each $X_k^n$ is an independent copy of $\det(Z(J_{nm}))$ with $n \geq m$ and $\rho > 0$. It holds that*

$$\mathrm{P}((1-\delta)\,\mathrm{per}\, J_{nm} \leq Y_n \leq (1+\delta)\,\mathrm{per}\, J_{nm}) \geq 1 - \frac{1}{\delta^2 n^\rho} \tag{4.2}$$

*for all $n \geq 2$.*

Both propositions are easily explained. The first is a consequence of the nice result of Silverstein [13] which says that if (4.1) holds, then $\lambda_1(Z(\tilde{J}_{nm_n}))$ converges in probability to a positive constant as $n \to \infty$. In other words, in this setting Condition 2.1 trivially holds for all $\varepsilon$ small enough.

Proposition 4.2 makes use of the well-known result (see again [13]) that the determinant of $Z(J_{nm})$ has the distribution $\chi_n^2 \chi_{n-1}^2 \cdots \chi_{n-m+1}^2$. Here the notation refers to the distribution of the product of independent random variables with the indicated $\chi^2$ distributions. A proof of this fact may be drawn from revisiting Lemma 3.1, as follows: Let $A$ denote an element of $\mathbf{M}(n, m, \mathbb{R}_+)$ and $A_k$ the matrix formed by removing the $k$th column. Again bring in the random matrix $X(A)$ with columns by $\bar{x}_1, \ldots, \bar{x}_m$ with $\bar{x}_k = D_k x_k$. Lemma (3.3) provides that

$$\det(Z(A)) = \det(Z(A_m))[x_m^T D_m P_m D_m x_m]$$
$$= \prod_{k=1}^m [x_k^T D_k P_k D_k x_k] = \prod_{k=1}^m \left[ \sum_{i=1}^{n-k+1} \lambda_{i+k-1}(D_k P_k D_k) \hat{x}_{ik}^2 \right], \quad (4.3)$$

in which the last equality is in law with the $\{\hat{x}_{ik}\}$ independent standard Gaussians. In the flat case $A = J_{nm}$ is affected by $D_k = I$ for all $k$, which is to say that $\lambda_i(P_k) = 1$ when $i \geq k$. The advertised distributional identity follows.

The use of this is in computing moments of $\det(J_{nm})$. That is,

$$\operatorname{per} J_{nm} = \mathrm{E}[\det Z(J_{nm})] = \mathrm{E}\left[ \prod_{k=n-m+1}^n \chi_k^2 \right] = \frac{n!}{(n-m)!},$$

which we knew before, but now also

$$\mathrm{E}[\det Z(J_{nm})]^2 = \mathrm{E}\left[ \prod_{k=n-m+1}^n (\chi_k^2)^2 \right]$$
$$= \prod_{k=n-m+1}^n (k^2 + 2k) = \left[ \frac{n!}{(n-m)!} \right]^2 \prod_{k=n-m+1}^n \left( 1 + \frac{2}{k} \right)$$
$$\leq \left[ \frac{n!}{(n-m)!} \right]^2 \exp\left( \sum_{i=k}^n \frac{2}{k} \right) \leq \frac{(n+1)^2}{4} \left[ \frac{n!}{(n-m)!} \right]^2.$$

With this estimate, the proof of Proposition 4.2 follows easily from Chebychev's inequality.

The question posed here is whether either approach [taking advantage of either the restrictive geometry as in (4.1) or the determinantal formula (4.3)] might lead to similarly sharp concentration in the more general $\mathbf{M}(n, m, [a, b])$ case. Believing that this should be so really comes down to believing that the

16         S. FRIEDLAND, B. RIDER AND O. ZEITOUNIbottom of the spectrum of $Z(A)$ for $A \in \mathbf{M}(n,m,[a,b])$ is not much worse than that of $Z(J_{nm})$. Of course, providing support for the latter statement has been the main technical goal of the present work.

## APPENDIX

A large part of the above argument entailed proving a certain integrability of the logarithm of the small eigenvalues of a Wishart-type matrix. It is interesting that the issue of controlling the bottom of the spectrum comes up in a great many problems (see, once more, [1] for an example). While not directly relevant for the study of the permanent, we wish to point out in this appendix that an exact analysis of the flat case $(J_n)$ reveals a much stronger integrability than that proved in Proposition 3.1. It is natural (and an underlying theme of this paper) to suppose the actuality of the more general case $\mathbf{M}(n,[a,b])$ is similar. For brevity we present the computation in the complex setting (the computation in the real case employs Pfaffians and requires nontrivial modifications), where $Y(\tilde{J}_n)_{ij} = (\frac{x_{ij}^R + \sqrt{-1} x_{ij}^I}{\sqrt{2n}})$ with $x_{ij}^{R,I}$ independent standard Gaussians. In fact, we note that in this case [12] have computed the law of the determinant of $Y(\tilde{J}_n)^* Y(\tilde{J}_n)$.

Our result is the following.

PROPOSITION A.1. *For $Y(\tilde{J}_n)$, an $n \times n$ matrix with entries independent complex Gaussians of mean 0 and variance $1/\sqrt{n}$, the eigenvalues $\lambda_i$ of $Y(\tilde{J}_n)^* Y(\tilde{J}_n)$ satisfy*

$$\lim_{\varepsilon \to 0} \lim_{n \to \infty} \mathrm{E}\left[\frac{1}{n} \sum_{\lambda_i < \varepsilon} \lambda_i^{-\alpha}\right] = 0$$

*for any $\alpha < 1/2$.*

PROOF. The present ensemble is integrable in the sense that the joint density of the eigenvalues $\lambda_1, \ldots, \lambda_n$ is known explicitly [9]:

$$\begin{aligned}
P(\lambda_1, \lambda_2, \ldots, \lambda_n) \\
&= C_n \exp\left[-n \sum_{i=1}^n \lambda_i\right] \prod_{i<j} (\lambda_i - \lambda_j)^2 \\
&= \frac{1}{n!} \det\left[e^{-n\lambda_i/2 - n\lambda_j/2} \sum_{k=0}^{n-1} \sqrt{n} L_k^0(n\lambda_i) \sqrt{n} L_k^0(n\lambda_j)\right]_{0 \leq i,j \leq n-1},
\end{aligned} \quad \text{(A.1)}$$

where $L_k^0$ denote the Laguerre polynomials: the family $\{L_k^\beta\}$ orthogonalized $x^{k+\beta} e^{-x}$ on $[0, \infty)$. From the determinantal formula (A.1), you may derive

the eigenvalue density

$$p_n(x) \equiv E\left[\frac{1}{n}\sum_{k=1}^{n}\delta_{\lambda_k}(x)\right] = e^{-nx}\sum_{k=0}^{n-1}(L_k^0(nx))^2.$$

By Christoffel–Darboux and the rules

$$\frac{d}{dx}L_k^\beta(x) = -L_{k-1}^{\beta+1} \quad \text{and} \quad L_{k-1}^\beta(x) = L_{k-1}^{\beta+1}(x) + L_{k-2}^{\beta+1}(x),$$

there is also the form

$$p_n(x) = ne^{-nx}\left[\frac{d}{dx}L_n^0(nx)L_{n-1}^0(nx) - L_n^0(nx)\frac{d}{dx}L_{n-1}^0(nx)\right]$$

$$= ne^{-nx}[L_n^1(nx)L_{n-2}^1(nx) - (L_{n-1}^1(nx))^2].$$

Thus, the integral to be examined is

(A.2) $$n\int_0^\varepsilon x^{-\alpha}[L_n^1(nx)L_{n-2}^1(nx) - (L_{n-1}^1(nx))^2]e^{-nx}\,dx.$$

Near zero, it is known [14] that $e^{-x/2}L_n^1(x) \leq Cnx^2$ for $0 < x < K/n$ with a large constant $K$, which allows you to dispel of the integral (A.2) for $x < K/n^2$: either term is of order

$$n\int_0^{K/n^2} x^{-\alpha}(e^{-nx/2}L_n^1(nx))^2\,dx \leq Cn^{2+\alpha}\int_0^{K/n} x^{4-\alpha}\,dx = O(n^{-3+2\alpha}).$$

For what remains, one needs the following (see [5]): uniformly on $0 < z < c\nu$ ($c < 1, \nu = 4n+4$),

(A.3)
$$e^{-z/2}L_n^1(z)$$
$$= (n-1)\frac{1}{z}\left(\frac{\psi}{\psi'}(z/n)\right)^{1/2}\left[J_1(\nu\psi(z/n)) + O\left(\frac{\sqrt{z}}{n^{3/2}}f(\nu\psi(z/n))\right)\right],$$

where $\psi(t) = (1/2)\sqrt{t-t^2}+(1/2)\sin^{-1}\sqrt{t}$, $f(t) = t$ for $t < 1$, $t^{-1/2}$ otherwise and $J_1$ is the Bessel function. We consider $z = nx, x < \varepsilon \ll 1$, on which $\psi'/\psi(z/n) \sim \sqrt{x}$. Note also that

$$\sqrt{\pi/2}J_1(z) \simeq \frac{1}{\sqrt{z}}\cos(z - 3\pi/4) \qquad \text{for } z\uparrow\infty.$$

Substituting (A.3) into the restricted integral (A.2), we first consider terms involving the second factor in (A.3). On the range $K/n^2 \leq x \leq \varepsilon$ we have $\frac{\sqrt{z}}{n^{3/2}}f(\nu\psi(z/n)) \leq cn^{-3/2}x^{1/4}$ and $J_1(\nu\psi(z/n)) \leq c/\sqrt{nx}$, yielding the contributions of order

$$n\int_{K/n^2}^\varepsilon x^{-\alpha}(1/\sqrt{x})^2(x^{1/4}/n^{3/2})^2\,dx = O(n^{-2+2\alpha})$$

and

$$n \int_{K/n^2}^{\varepsilon} x^{-\alpha}(1/\sqrt{x}\,)^2(1/n^{1/2}x^{1/4})(x^{1/4}/n^{3/2})\,dx = O(n^{-1+2\alpha}),$$

both vanishing for $n \to \infty$ as soon as $\alpha < 1/2$. That leaves us with

$$\int_0^{\varepsilon} x^{-\alpha} p_n(x)\,dx$$
$$= n \int_{K/n^2}^{\varepsilon} x^{-\alpha}\frac{1}{x}[J_1((n+1)\sqrt{x}\,)J_1((n-1)\sqrt{x}\,) - J_1^2(n\sqrt{x}\,)]\,dx$$
$$- \frac{1}{n}\int_{K/n^2}^{\varepsilon} x^{-\alpha}\frac{1}{x}J_1((n+1)\sqrt{x}\,)J_1((n-1)\sqrt{x}\,)\,dx + O(n^{-1+2\alpha}).$$

For the first term on the right-hand side, the integrand is overestimated as in $|\cos((n+1)\sqrt{x}\,)\cos((n-1)\sqrt{x}\,) - \cos^2(n\sqrt{x}\,)| \leq x$ for $0 \leq x \leq \varepsilon$. The remaining integral is then controlled by a constant multiple of $\int_{K/n^2}^{\varepsilon} x^{-1/2-\alpha}\,dx = \varepsilon^{1/2-\alpha} - O(n^{2\alpha-1})$. The second term is even easier: the bound $J_1(nz) \leq C/\sqrt{nz}$ shows it to be of order $(1/n^2)\int_{K/n^2}^{\varepsilon} x^{-3/2-\alpha}\,dx \sim n^{2\alpha-1}$. $\square$

S. FRIEDLAND
DEPARTMENT OF MATHEMATICS
  STATISTICS AND COMPUTER SCIENCE
UNIVERSITY OF ILLINOIS AT CHICAGO
851 S. MORGAN ST.
CHICAGO, ILLINOIS 60607-7045
USA
E-MAIL: friedlan@uic.edu
URL: www.math.uic.edu/~friedlan

B. RIDER
DEPARTMENT OF MATHEMATICS
DUKE UNIVERSITY
BOX 90320
DURHAM, NORTH CAROLINA 27708-0320
USA
E-MAIL: rider@math.duke.edu
URL: www.math.duke.edu/~rider

O. ZEITOUNI
DEPARTMENTS OF MATHEMATICS
  AND OF ELECTRICAL ENGINEERING
TECHNION, HAIFA 32000
ISRAEL
AND
SCHOOL OF MATHEMATICS
UNIVERSITY OF MINNESOTA
MINNEAPOLIS, MINNESOTA 55455
USA
E-MAIL: zeitouni@math.umn.edu
URL: www-ee.technion.ac.il/~zeitouni